\documentclass[10pt,twoside]{article}

\usepackage[centertags]{amsmath}
\usepackage{amsfonts}
\usepackage{amssymb}
\usepackage{amsthm}
\usepackage{epsfig}
\usepackage{color}
\usepackage[titletoc,toc,title]{appendix}

\allowdisplaybreaks[4]

\allowdisplaybreaks[4]

\definecolor{c20}{rgb}{0.,0.7,0.}
\definecolor{c30}{rgb}{0.,0.,1.}
\definecolor{c40}{rgb}{1,0.1,0.7}
\definecolor{c50}{rgb}{1,0,0}
\definecolor{c60}{rgb}{1,0.9,0.1}

\def\raini{\to {\boldsymbol{\infty}}}

\def\indi{{1\kern-.20em\rm I}}

\newtheorem{theorem}{Theorem}[section]
\newtheorem{definition}{Definition}[section]

\newtheorem{lem}{Lemma}[section]
\newtheorem{prop}{Proposition}[section]
\newtheorem{corollary}{Corollary}[section]
\newtheorem{rem}{Remark}[section]
\newtheorem{ex}{Example}[section]

\begin{document}
\baselineskip 15pt \setcounter{page}{1}
\title{\bf \Large  Almost sure convergence for the maximum of nonstationary random fields}
\author{{\small Lu\'{i}sa Pereira$^{1}$\footnote{ lpereira@ubi.pt}\ \ and Zhongquan Tan$^{2}$}\\
\\
{\small\it  1. Department of Mathematics, University of Beira Interior, Portugal}\\
{\small\it 2. College of Mathematics, Physics and Information
Engineering, Jiaxing University, Jiaxing 314001, PR China}\\
}
 \maketitle
 \baselineskip 15pt

\begin{quote}
{\bf Abstract:}\ \ We obtain an almost sure limit theorem for the maximum of nonstationary random fields under some dependence conditions.
The obtained result is applied to Gaussian random fields.



\end{quote}

\section{Introduction}

In recent years various authors discussed almost sure versions of distributional limit theorems. The first result on Almost Sure Central Limit Theorem (ASCLT) presented independently by Brosamler (1988), Schatte (1988) and Lacey and Philipp (1990) extended the classical central limit theorem to an almost sure version.

For an i.i.d. sequence $\{X_n\}_{n\in \mathbb N}$ with zero mean, unit variance and partial sum $S_k=\sum_{i=1}^k X_i$, $k\geq 1$, the simplest version of the ASCLT states that
$$
\frac{1}{\log n}\sum_{k=1}^n\frac{1}{k}\indi_{\left\{S_k\leq \sqrt{k}x\right\}}\rightarrow\Phi(x) \ \ a.s.,
$$
for any fixed $x \in \mathbb{R}$, where $a.s.$ means almost surely, $\indi_A$ denotes the indicator function of the event $A$ and $\Phi(x)$ is the standard normal distribution function.

Later on the ASCLTs for some other functions of random variables were studied. Namely, in  Fahrner and Stadm\"{u}ller (1998), Cheng, Peng and Qi (1998) and Berkes and Cs\'{a}ki (2001) the ASCLTs for the maximum of an i.i.d. random sequence were proved.

Let $\{X_n\}_{n\in \mathbb N}$ be an i.i.d. sequence, and let $M_k=\max_{1\leq i \leq k }X_i$ denotes the partial maximum, $k\geq 1$. If there exist normalizing constants $a_k>0$, $b_k \in \mathbb{R}$ and a nondegenerate distribution function $G(x)$ such that
$$
P\left(M_n\leq a_nx+b_n\right)\rightarrow G(x),
$$
then we have
$$
\frac{1}{\log n}\sum_{k=1}^n \frac{1}{k} \indi_{\left\{M_k\leq a_kx+b_k\right\}}\rightarrow G(x) \ \ a.s., 
$$
for any continuity point $x$ of $G$.
It is well known that $G(x)$ must be of the same type as the extreme value distribution $G(x)=\exp\left\{-(1+\gamma x)^{-\frac{1}{\gamma}}\right\}$, where $\gamma$ is the so-called extreme value index.

On the other hand, the ASCLTs for the maximum of some dependent, stationary normal sequences were obtained by Cs\'{a}ki and Gonchigdanzan (2002), while the ASCLT for the maximum of some dependent, but not necessarily stationary sequences was established by  Peng and Nadarajah (2011) and Chen and Lin (2006).\\
So far many results have been obtained for the a.s. convergence for extremes of random sequences but few for random fields. Some works which are worthwhile to mention in this place are the papers of Choi (2010) and Tan and Wang (2014), where they established the ASCLT for the maximum of stationary and nonstationary normal random fields, respectively. Random fields are of increasing interest in applications such as environmental assessment over entire regions of space.

In this paper we prove an ASCLT for the maximum of nonstationary random fields, $\mathbf{X}=\left\{ X_{\mathbf{n}}\right\} _{{\mathbf n}\in \mathbb{Z}_{+}^2}$, where $\mathbb{Z}_{+}$ is the set of all positive integers,subject to conditions on long range and local dependencies.

Throughout we shall say that the pair $(\textbf{I},\textbf{J})$, $\textbf{I},\textbf{J}\subseteq \mathbb{Z}_{+}^{2}$, is in $\mathcal{S}_{i}(l)$, for each $i=1,2$, if the distance between $\Pi _{i}(\mathbf{I})$ and
$\Pi _{i}(\mathbf{J})$ is greater or equal to $l$, where $\Pi _{i},\,i=1,2$,
denote the cartesian projections. For $\mathbf{i}=(i_{1},i_{2})$ and $\mathbf{j}%
=(j_{1},j_{2}),$ $\mathbf{i}\leq \mathbf{j}$ means $i_{k}\leq j_{k},$ $k=1,2$, and $\mathbf{n}=(n_1, n_2)\rightarrow\boldsymbol{\infty}$ means $n_k\rightarrow \infty, k=1,2$.
Considering that $\{u_{\mathbf{n},\mathbf{i}}:\mathbf{i}\leq \mathbf{n}\}_{\mathbf{n}\geq
\mathbf{1}}$ is a sequence of real numbers and $\mathbf{I}$ a subset of the rectangle of points $\mathbf{R}%
_{\mathbf{n}}=\{1,\ldots ,n_{1}\}\times \{1,\ldots ,n_{2}\}$, we will denote
the event $\{X_{\mathbf{i}}\leq u_{\mathbf{n},\mathbf{i}}:\mathbf{i}\in \mathbf{I}\}$
by $\left\{M_{\mathbf{n}}(\mathbf{I})\leq u\right\}$ or simply by $\left\{M_{\mathbf{n}}\leq u\right\}$ when $\mathbf{I%
}=\mathbf{R}_{\mathbf{n}}$. Let $\mathbf{1}=(1,1)$.

As discussed in Pereira and Ferreira (2005,2006) in order to prove that the probability of no exceedances of high values over $\mathbf{R}%
_{\mathbf{n}}$ can be approximated by $\exp\{-\tau\}$, where $\tau$ is the limiting mean number of exceedances, the following conditions are needed.

The first is a coordinatewise-mixing type condition as the $\Delta (u_{\mathbf{n}})-$ condition introduced in Leadbetter and
Rootz\'{e}n (1998), which restrict dependence by limiting $$\left|P\left(M_{\bf n}({\bf I}_1)\leq u,M_{\bf n}({\bf I}_2)\leq u\right)-P\left(M_{\bf n}({\bf I}_1)\leq u\right)P\left(M_{\bf n}({\bf I}_1)\leq u\right)\right|$$ with the two indexes sets ${\bf I}_1$ and ${\bf I}_2$ being "separated" from each other by a certain distance along each coordinate direction.

\begin{definition}
Let $\mathcal{F}$ be a family of indexes sets in $\mathbf{R}_{\mathbf{n}}.$
The nonstationary random field $\mathbf{X}$ on $\mathbb{Z}_{+}^{2}$
satisfies the condition $D\mathbb{(}u_{\mathbf{n},\mathbf{i}}\mathbb{)}$
over $\mathcal{F}$ if there exist sequences of integer valued constants $%
\left\{ k_{n_{i}}\right\} _{n_{i}\geq 1},$\textit{\ }$\left\{
l_{n_{i}}\right\} _{n_{i}\geq 1},$\textit{\ }$i=1,2,$ such that, as $\mathbf{%
n}=({n_{1},n_{2})}\longrightarrow \mathbf{\boldsymbol{\infty} ,}$ we have
\begin{equation}
\left( k_{n_{1}},k_{n_{2}}\right) \longrightarrow \mathbf{\boldsymbol{\infty} ,}\text{ }%
\left( \tfrac{k_{n_{1}}l_{n_{1}}}{n_{1}},\tfrac{k_{n_{2}}l_{n_{2}}}{n_{2}}%
\right) \longrightarrow \mathbf{0}  \notag
\end{equation}
and $\left( k_{n_{1}}\Delta _{{\bf{n}},l_{n_{1}}}^{(1)},k_{n_{1}}k_{n_{2}}\Delta
_{{\bf{n}},l_{n_{2}}}^{(2)}\right) \longrightarrow \mathbf{0}$, where $\Delta
_{{\bf{n}},l_{n_{i}}}^{(i)},i=1,2,$ are the components of the mixing coefficient,
defined as follows: {\small {\
\begin{equation*}
\Delta _{{\bf{n}},l_{n_{1}}}^{(1)}=\sup \left| P\left( M_{\mathbf{n}}(\mathbf{I}%
_{1})\leq u,M_{\mathbf{n}}(\mathbf{I}_{2})\leq u%
\right) -P\left( M_{\mathbf{n}}(\mathbf{I}_{1})\leq u\right)
P\left( M_{\mathbf{n}}(\mathbf{I}_{2})\leq u\right) \right| ,
\end{equation*}
}}where the supremum is taken over pairs of $\mathbf{I}_{1}$ and $\mathbf{I}%
_{2}$ in $S_{1}(l_{n_{1}})\cap \mathcal{F}$ $,${\small {\
\begin{equation*}
\Delta _{{\bf n},l_{n_{2}}}^{(2)}=\sup \left| P\left( M_{\mathbf{n}}(\mathbf{I}%
_{1})\leq u,M_{\mathbf{n}}(\mathbf{I}_{2})\leq u%
\right) -P\left( M_{\mathbf{n}}(\mathbf{I}_{1})\leq u\right)
P\left( M_{\mathbf{n}}(\mathbf{I}_{2})\leq u\right) \right| ,
\end{equation*}
}}where the supremum is taken over pairs of $\mathbf{I}_{1}$ and $\mathbf{I}%
_{2}$ in $S_{2}(l_{n_{2}})\cap \mathcal{F}.$
\end{definition}

This condition was used to guarantee the asymptotic independence for maxima over disjoint rectangles of indexes (Pereira and Ferreira (2006)) which is a fundamental result for extending some results of the extreme value theory of stationary random fields to nonstationary case.

\begin{prop}  Suppose that the random field $\mathbf{X}$ satisfies the condition  $D\mathbb{(}u_{\mathbf{n},\mathbf{i}}\mathbb{)}$ over $\mathcal{F}$ such that $(\mathbf{I}\subset\mathbf{J}\wedge\mathbf{J}\in\mathcal{F})\Rightarrow \mathbf{I}\in\mathcal{F}$ and for $\{u_{\mathbf{n},\mathbf{i}}:\mathbf{i}\leq \mathbf{n}\}_{\mathbf{n}\geq
\mathbf{1}}$ such that
\begin{equation*}
\left\{n_1n_2\max\left\{P(X_{\mathbf{i}}>u_{\mathbf{n},\mathbf{i}}):\mathbf{i}\leq\mathbf{n}\right\}\right\}_{\mathbf{n}\geq\mathbf{1}}\ is\ \ bounded.
\end{equation*}
If $\mathbf{V}_{r,p}=I_r\times J_{r,p},\ r=1,\ldots,k_{n_1},\ p=1,\ldots,k_{n_2},$ are disjoint rectangles in $\mathcal{F}$, then, as $\mathbf{n}\rightarrow \boldsymbol{\infty}$,
\begin{equation*}
P\left(\underset{r,p}{\bigcap }\left\{M_{\mathbf {n}}(\mathbf{V}_{r,p})\leq u\right\}\right)-\underset{r,p}{\prod}P(M_{\mathbf {n}}(\mathbf{V}_{r,p})\leq u)\rightarrow 0.
\end{equation*}
\end{prop}

In Pereira and Ferreira (2005), in addition to the coordinatewise-mixing condition, it is restricted the local path behaviour with respect to exceedances. It is used the idea of Leadbetter and Rootz\'{e}n (1998) in combination with H\"{u}sler (1986) to generalize to the nonstationary case a local dependence condition, $D^{\prime }(u_{\mathbf{n},\mathbf{i}})$, that avoids clustering of exceedances of $u_{\mathbf{n},\mathbf{i}}$.

\begin{definition}
Let $\mathcal{E}(u_{\mathbf{n},\mathbf{i}})$ denote the family of
indexes sets $\mathbf{I}$ such that
\begin{equation*}
\underset{\mathbf{i}\in \mathbf{I}}{\sum }P\left( X_{\mathbf{i}}>u_{\mathbf{n%
},\mathbf{i}}\right) \leq \frac{1}{k_{n_{1}}k_{n_{2}}}\underset{\mathbf{i}%
\leq \mathbf{n}}{\sum }P\left( X_{\mathbf{i}}>u_{\mathbf{n},\mathbf{i}%
}\right) .
\end{equation*}
\textit{The} \textit{condition} $D^{\prime }(u_{\mathbf{n},\mathbf{i}})$
holds for $\mathbf{X}$ \textit{if, for each }$\mathbf{I}\in \mathcal{E}(u_{%
\mathbf{n},\mathbf{i}})$, we have, as  $\mathbf{n}\rightarrow \boldsymbol{\infty ,}$ \textit{\ }
\begin{equation*}
k_{n_{1}}k_{n_{2}}\underset{\mathbf{i,j}\in \mathbf{I}}{\sum }P(X_{\mathbf{i}%
}>u_{\mathbf{n},\mathbf{i}},X_{\mathbf{j}}>u_{\mathbf{n},\mathbf{j}})\longrightarrow 0.
\end{equation*}
\end{definition}

That condition, which bounds the probability of more than one exceedance above the levels $u_{\mathbf{n},\mathbf{i}}$ in a rectangle with a few indexes, and the coordinatewise-mixing $D(u_{{\bf n},{\bf i}})$ condition lead to a Poisson approximation for the probability of no exceedances over $\mathbf{R}_{%
\mathbf{n}}$ (see, Pereira and Ferreira (2005)).

\begin{prop}
\textit{Suppose that the} nonstationary\textit{\ random field} $\mathbf{X}$
\textit{satisfies} $D(u_{\mathbf{n},\mathbf{i}})$ and $D^{\prime }(u_{\mathbf{n},\mathbf{i}})$ over $\mathcal{E}(u_{\mathbf{n},\mathbf{i}})$ and
\begin{equation*}
\left\{ n_{1}n_{2}\max\left\{ P\left( X_{\mathbf{i}}>u_{\mathbf{n},\mathbf{i}%
}\right) :\mathbf{i}\leq \mathbf{n}\right\} \right\} _{\mathbf{n}\geq
\mathbf{1}}\text{ \ \ \ \ is bounded.}
\end{equation*}
Then,
\begin{equation}
P(M_{\mathbf{n}}\leq u_{\mathbf{n},\mathbf{i}})\xrightarrow
[\mathbf{n}\raini]{}\exp (-\tau ),\text{ \ \ }\tau >0,  \notag
\end{equation}
\textit{if and only if}
\begin{equation}
\underset{\mathbf{i}\leq \mathbf{n}}{\sum }P(X_{\mathbf{i}}>u_{\mathbf{n},%
\mathbf{i}})\xrightarrow [\mathbf{n}\raini]{}\tau .  \notag
\end{equation}
\end{prop}
\vspace{1cm}

The a.s. version of Proposition 1.2 is given in Section 2. Section 3 is devoted to the a.s. convergence for the maximum of a normal random field. We prove that our main results are more general than the results established in Choi (2010) and Tan and Wang (2014). All the proofs are collected in appendices.

\section{Main result}

Throughout the paper $<\!<$ stands for $a=O(b)$.

In order to formulate the main result we need to strengthen condition $D\mathbb{(}u_{\mathbf{n},\mathbf{i}}\mathbb{)}$ as follows.

\begin{definition}
Let $\mathcal{F}$ be a family of indexes sets in $\mathbf{R}_{\mathbf{n}}.$
The nonstationary random field $\mathbf{X}$ on $\mathbb{Z}_{+}^{2}$
satisfies the condition $D^*\mathbb{(}u_{\mathbf{n},\mathbf{i}}\mathbb{)}$
over $\mathcal{F}$ if there exist sequences of integer valued constants $%
\left\{ k_{n_{i}}\right\} _{n_{i}\geq 1},$\textit{\ }$\left\{
m_{n_{i}}\right\} _{n_{i}\geq 1},$\textit{\ }$i=1,2,$ such that, as $\mathbf{%
n}=({n_{1},n_{2})}\rightarrow \mathbf{\boldsymbol{\infty} ,}$ we have
\begin{equation}
\left( k_{n_{1}},k_{n_{2}}\right) \longrightarrow \mathbf{\boldsymbol{\infty} ,}\text{ }%
\left( \tfrac{k_{n_{1}}m_{n_{1}}}{n_{1}},\tfrac{k_{n_{2}}m_{n_{2}}}{n_{2}}%
\right) \longrightarrow \mathbf{0}  \notag
\end{equation}
and for some $\epsilon>0$
$$\alpha_{{\bf{n}},m_{n_{1}},m_{n_{2}}}=\sup_{\mathbf{1}\leq \mathbf{k}\leq \mathbf{n}}\alpha_{{\bf{n,k}},m_{n_{1}},m_{n_{2}}}<\!<(\log n_1 \log n_2)^{-(1+\epsilon)},$$  where $\alpha_{{\bf{n,k,l}},m_{n_{1}},m_{n_{2}}}$ is the mixing coefficient, defined as follows: {\small {\
\begin{equation*}
\alpha_{{\bf{n,k}},m_{n_{1}},m_{n_{2}}}=\sup_{({\bf I,  J})\in \mathcal{S}(m_{n_1}, m_{n_2})} \left| P\left( \bigcap_{{\bf i}\in {\bf I}}\left\{X_{\bf i}\leq u_{{\bf k},{\bf i}}\right\},\bigcap_{{\bf j}\in {\bf J}}\left\{X_{\bf j}\leq u_{{\bf n},{\bf j}}\right\}
\right) -P\left( \bigcap_{{\bf i}\in {\bf I}}\left\{X_{\bf i}\leq u_{{\bf k},{\bf i}}\right\}
\right)
P\left( \bigcap_{{\bf j}\in {\bf J}}\left\{X_{\bf j}\leq u_{{\bf n},{\bf j}}\right\}
\right) \right|,
\end{equation*}
}}
$\mathbf{k}=S(\mathbf{I})$ and
\begin{eqnarray*}
&&\mathcal{S}(m_{n_1},m_{n_2})=\left\{({\bf I},{\bf J})\subseteq \mathbf{R_{n}}^2:s(\Pi_2({\bf J}))-S(\Pi_2({\bf I}))\geq m_{n_2}
\vee s(\Pi_1({\bf J}))-S(\Pi_1({\bf I}))\geq m_{n_1}\right\},
\end{eqnarray*}
  with $S({\bf I})=\sup\left\{{\bf i}:{\bf i}\in {\bf I}\right\}$ and $s({\bf I})=\inf\left\{{\bf i}:{\bf i}\in {\bf I}\right\}$.
\end{definition}

\begin{theorem} Let ${\bf X}$ be a nonstationary random field satisfying conditions $D^*(u_{\mathbf{n},\mathbf{i}})$ and $D^{\prime }(u_{\mathbf{n},\mathbf{i}})$ over $\varepsilon (u_{{\bf n}, {\bf i}})$.
Assume that
$$
\underset{\mathbf{i}\leq \mathbf{n}}{\sum }P(X_{\mathbf{i}}>u_{\mathbf{n},%
\mathbf{i}})\xrightarrow [\mathbf{n}\raini]{}\tau, \ \ {\text for\ \ some \ \  } 0\leq\tau<\infty,
$$
and $\left\{ n_{1}n_{2}\max\left\{ P\left( X_{\mathbf{i}}>u_{\mathbf{n},\mathbf{i}%
}\right) :\mathbf{i}\leq \mathbf{n}\right\} \right\} _{\mathbf{n}\geq
\mathbf{1}}$ is bounded.
Then
\begin{equation}
\lim_{{\bf n}\rightarrow {\boldsymbol{\infty}}}\frac{1}{\log n_1 \log n_2}\sum_{{\bf k}\in {\bf R_n}}\frac{1}{k_1k_2}\indi_{\left\{\bigcap_{{\bf i}\leq {\bf k}}\left\{X_{\bf i}\leq u_{{\bf k},{\bf i}}\right\}\right\}}=\exp(-\tau) \ \ a.s.
\end{equation}
\end{theorem}

For stationary random fields, based on condition $D^{\prime }(u_{\mathbf{n}})$ in Leadbetter and Rootz\'{e}n (1998) and condition $D^*\mathbb{(}u_{\mathbf{n},\mathbf{i}}\mathbb{)}$ with $u_{{\bf n},{\bf i}}=u_{\bf n}$ we have the following result.

\begin{corollary}
Let ${\bf X}$ be a stationary random field satisfying conditions $D^{\prime }(u_{\mathbf{n}})$ and $D^*(u_{\mathbf{n}})$.
If
$$
n_1n_2P(X_{\mathbf{1}}>u_{\mathbf{n}})\xrightarrow [\mathbf{n}\raini]{}\tau, \ \ {\text for\ \ some \ \ } 0\leq\tau<\infty,
$$
then
$$
\lim_{{\bf n}\rightarrow {\boldsymbol{\infty}}}\frac{1}{\log n_1 \log n_2}\sum_{{\bf k}\in {\bf R_n}}\frac{1}{k_1k_2}\indi_{\left\{\bigcap_{{\bf i}\leq {\bf k}}\left\{X_{\bf i}\leq u_{{\bf k}}\right\}\right\}}=\exp(-\tau) \ \ a.s.
$$
\end{corollary}

\section{Normal random fields}



Normality occupies a central place in probability and statistical theory, and a most important class of random fields consists of those which are normal. Their importance is enhanced by the fact that the specification of their finite-dimensional distributions are simple, they are reasonable models for many natural phenomenon, estimation and inference are simple and the model is specified by expectations and covariances.

The almost sure convergence for the maximum of a normal random field is investigated. The covariance conditions given by Tan and Wang (2014) for the a.s. convergence given in (1) is compared with the dependence conditions used in Section 2. An example satisfying the conditions of Theorem 2.1 but not
the conditions in Tan and Wang (2014) will be given.

\vspace{0.3cm}
Tan and Wang (2014) gave simple conditions on the covariances of nonstationary standardized normal random field to ensure that (1) holds.

\begin{theorem}
Let $\mathbf{X}$ be a non-stationary
standardized normal random field. Assume that the covariance functions $r_{\mathbf{i,j}}$ satisfy
$|r_{\mathbf{i,j}}|<\rho_{\mathbf{|i-j|}}$
for some sequence $\{\rho_{\mathbf{n}}\}_{\mathbf{n}\in \mathbb{N}^{2}-\{\mathbf{0}\}}$ such that for some $\epsilon>0$,
\begin{eqnarray}
\rho_{(n_{1},0)}<\!<(\log n_{1})^{-(1+\epsilon)},\ \
\rho_{(0,n_{2})}<\!<(\log n_{2})^{-(1+\epsilon)}, \ \
\rho_{\mathbf{n}}<\!<(\log n_{1}n_{2})^{-(1+\epsilon)},
\end{eqnarray}
and $\sup_{\mathbf{n}\in \mathbb{N}^{2}-\{\mathbf{0}\}}|\rho_{\mathbf{n}}|<1$ hold.
Let the constants $\{u_{\mathbf{n,i}},\mathbf{i}\leq \mathbf{n}\}_{\mathbf{n}\geq \mathbf{1}}$ be such that
$n_{1}n_{2}(1-\Phi(\lambda_{\mathbf{n}}))$ is bounded, where
$\lambda_{\mathbf{n}}=\min_{\mathbf{i}\in \mathbf{R_{n}}}u_{\mathbf{n,i}}$.
Suppose that $\lim_{\mathbf{n}\rightarrow\infty}\sum_{\mathbf{i}\in \mathbf{R_{n}}}(1-\Phi(u_{\mathbf{n,i}}))=\tau\in[0,\infty)$.
Then, the assertion of Theorem 2.1 holds.
\end{theorem}
\vspace{0.2cm}

So, we have the a.s. convergence given in (1) under $D^*(u_{\mathbf{n},\mathbf{i}})$ and $D^{\prime }(u_{\mathbf{n},\mathbf{i}})$ conditions as under the covariance conditions given in (2). Hence it is desirable to investigate the relation of condition (2) and the conditions $D^*(u_{\mathbf{n},\mathbf{i}})$ and $D^{\prime }(u_{\mathbf{n},\mathbf{i}})$. This relation will be described through Theorem 3.2 and an example.

\vspace{0.2cm}

\begin{theorem}
Let $\mathbf{X}$ be a nonstationary
standardized normal random field. Assume that the covariance functions $r_{\mathbf{i,j}}$ satisfy
$|r_{\mathbf{i,j}}|<\rho_{\mathbf{|i-j|}}$
for some sequence $\{\rho_{\mathbf{n}}\}_{\mathbf{n}\in \mathbb{N}^{2}-\{\mathbf{0}\}}$ verifying (2) and $\sup_{\mathbf{n}\in \mathbb{N}^{2}-\{\mathbf{0}\}}|\rho_{\mathbf{n}}|<1$.
Let the constants $\{u_{\mathbf{n,i}},\mathbf{i}\leq \mathbf{n}\}_{\mathbf{n}\geq \mathbf{1}}$ be such that
$n_{1}n_{2}(1-\Phi(\lambda_{\mathbf{n}}))$ is bounded, where
$\lambda_{\mathbf{n}}=\min_{\mathbf{i}\in \mathbf{R_{n}}}u_{\mathbf{n,i}}$.
Suppose that $\lim_{\mathbf{n}\rightarrow\infty}\sum_{\mathbf{i}\in \mathbf{R_{n}}}(1-\Phi(u_{\mathbf{n,i}}))=\tau\in[0,\infty)$ holds.
Then, ${\bf X}$ satisfies the conditions conditions $D^*(u_{\mathbf{n},\mathbf{i}})$ and $D^{\prime }(u_{\mathbf{n},\mathbf{i}})$ over $\varepsilon (u_{{\bf n}, {\bf i}})$.
\end{theorem}

Therefore, the result given in Theorem 2.1 for the particular case in which ${\bf X}$ is a normal random field, is a more general result than Theorem 3.1.

\begin{rem}
 The assertion of Theorem 3.2 still holds for stationary normal random fields with similar conditions on the correlation functions and $u_{\mathbf{n,i}}=u_{\mathbf{n}}$.

\end{rem}

Next, we give an example which satisfies conditions of Theorem 2.1 but not conditions of Theorem 3.1.
\begin{ex}
Let $X_{\mathbf{n}}$ be a stationary normal field with covariance function
$$\gamma_{\mathbf{n}}=\gamma_{(n_{1},n_{2})}=\prod_{i=1}^{2}\left(\left(1-\frac{|n_{i}|}{2e^{2}}\right)^{1/2}\omega(n_{i})\indi_{\{|n_{i}|\leq e^{2}\}}+\left(\frac{1}{\log n_{i}}\right)^{1/2}\omega(n_{i})\indi_{\{|n_{i}|>e^{2}\}}\right),$$
where $\omega(n)=\Pi_{j=1}^{\infty}\cos(\frac{1}{3^{j}}n)$. Let the constants $\{u_{\mathbf{n}}\}_{\mathbf{n}\geq \mathbf{1}}$ be such that
$n_{1}n_{2}(1-\Phi(u_{\mathbf{n}}))\rightarrow\tau\in[0,\infty)$.
Then $X_{\mathbf{n}}$ satisfies the conditions of Theorem 2.1  but not the conditions of Theorem 3.1.
\end{ex}

Choi (2002) has showed that
$$\gamma_{n}:=\left(1-\frac{|n|}{2e^{2}}\right)^{1/2}\omega(n)\indi_{\{|n|\leq e^{2}\}}+\left(\frac{1}{\log n}\right)^{1/2}\omega(n)\indi_{\{|n|>e^{2}\}}$$ is a covariance function and $$\limsup_{n\rightarrow\infty}\gamma_{n}\log n=\infty.$$
It is easy to see that $\gamma_{\mathbf{n}}=\gamma_{n_{1}}\gamma_{n_{2}}$ is a covariance function and
$$\limsup_{n_{1}\rightarrow\infty}\gamma_{(n_{1},0)}\log n_{1}=\infty,\ \ \limsup_{n_{2}\rightarrow\infty}\gamma_{(0,n_{2})}\log n_{2}=\infty$$
and
$$\limsup_{\mathbf{n}\rightarrow\infty}\gamma_{\mathbf{n}}\log (n_{1}n_{2})>0.$$
So, $\gamma_{\mathbf{n}}$ does not satisfies the conditions of Theorem 3.1. In Appendix B, we show $\gamma_{\mathbf{n}}$ satisfies the conditions of Theorem 2.1.




\begin{appendix}\setcounter{equation}{0}
\section*{Appendix A: Proofs for Section 2}\setcounter{section}{0}

\refstepcounter{section}


Let $B_{\bf k}({\bf R_k})=\bigcap_{\bf i\in R_k}\left\{X_{\bf i}\leq u_{\bf k,i}\right\}$ and $\overline{B}_{\bf k}({\bf R_k})=\bigcup_{\bf i\in R_k}\left\{X_{\bf i}> u_{\bf k,i}\right\}$. For ${\bf k,l\in R_n}$ such that ${\bf k}\neq {\bf l}$ and $u_{{\bf l},{\bf i}}\geq u_{{\bf k},{\bf i}}$, let $m_{l_{i}}=\log l_{i}$. Note that $k_{1}k_{2}\leq l_{1}l_{2}$.
Let $\mathbf{M^{*}}=\mathbf{M^{*}}_{\mathbf{kl}}=\mathbf{R_{k}}\cap \mathbf{R_{l}}$ and
$\mathbf{M}_{\mathbf{kl}}=\{(x_{1},x_{2}): (x_{1},x_{2})\in \mathbf{N}^{2}, 0\leq x_{i}\leq \sharp(\prod_{i}(\mathbf{M}^{*}))+m_{l_{i}}, i=1,2\}$,
where $\sharp$ denotes cardinality. Note that $\mathbf{M^{*}}\subset \mathbf{M_{kl}}$.

The proof of Theorem 2.1 will be given by means of several lemmas.

\begin{lem}
Let ${\bf X}$ be a nonstationary random field satisfying condition $D^*(u_{\mathbf{n},\mathbf{i}})$ over $\mathcal{F}$.
Assume that $\left\{ n_{1}n_{2}\max\left\{ P\left( X_{\mathbf{i}}>u_{\mathbf{n},\mathbf{i}%
}\right) :\mathbf{i}\leq \mathbf{n}\right\} \right\} _{\mathbf{n}\geq
\mathbf{1}}$ is bounded and $\alpha_{\mathbf{l},m_{l_1},m_{l_2}}<\!<(\log l_1 \log l_2)^{-(\epsilon+1)}$. Then, for ${\bf k,l\in R_n}$ such that ${\bf k}\neq {\bf l}$ and $u_{{\bf l},{\bf i}}\geq u_{{\bf k},{\bf i}}$ 
$$
\left|Cov(\indi_{\left\{\bigcap_{\bf i\in R_k}\left\{X_{\bf i}\leq u_{\bf k,i}\right\}\right\}}, \indi_{\left\{\bigcap_{\bf i\in R_l-R_k}\left\{X_{\bf i}\leq u_{\bf l,i}\right\}\right\}}) \right|<\!<\alpha_{\mathbf{l,k},m_{l_1},m_{l_2}}+\frac{m_{l_1}k_{2}}{l_1l_{2}}+\frac{m_{l_2}k_{1}}{l_{1}l_2}.
$$
\end{lem}
\textbf{Proof:} Write
\begin{eqnarray*}
&&\left|Cov\left(\indi_{\left\{\bigcap_{\bf i\in R_k}\left\{X_{\bf i}\leq u_{\bf k,i}\right\}\right\}}, \indi_{\left\{\bigcap_{\bf i\in R_l-R_k}\left\{X_{\bf i}\leq u_{\bf l,i}\right\}\right\}}\right) \right|\\
&=&\left|P(B_{\bf k}({\bf R_k})\cap B_{\bf l}({\bf R_l-R_k}))-P(B_{\bf k}({\bf R_k}))P(B_{\bf l}(\bf R_l-R_k))\right|\\
&\leq&\left|P(B_{\bf k}({\bf R_k})\cap B_{\bf l}({\bf R_l-R_k}))-P(B_{{\bf k}}({\bf R_k})\cap B_{\bf l}(\bf R_l-M_{kl}))\right|\\
&&+\left|P(B_{\bf k}({\bf R_k})\cap B_{\bf l}({\bf R_l-M_{kl}}))-P(B_{{\bf k}}({\bf R_k}))P(B_{\bf l}(\bf R_l-M_{kl}))\right|\\
&&+\left|P(B_{{\bf k}}({\bf R_k}))P(B_{\bf l}({\bf R_l-M_{kl}}))-P(B_{{\bf k}}({\bf R_k}))P(B_{\bf l}(\bf R_l-R_{k}))\right|\\
&=:&I_1+I_2+I_{3}.
\end{eqnarray*}
Using the condition that $\left\{ n_{1}n_{2}\max\left\{ P\left( X_{\mathbf{i}}>u_{\mathbf{n},\mathbf{i}%
}\right) :\mathbf{i}\leq \mathbf{n}\right\} \right\} _{\mathbf{n}\geq
\mathbf{1}}$ is bounded we get
\begin{eqnarray*}
I_1&=&\left|P(B_{\bf k}({\bf R_k})\cap B_{\bf l}({\bf R_l-R_k}))-P(B_{\bf k}({\bf R_k})\cap B_{{\bf l}}(\bf R_l-M_{kl}))\right|\\
&\leq& \left|P(B_{\bf l}({\bf R_l-R_k}))-P(B_{\bf l}(\bf R_l-M_{kl}))\right|\\
&\leq& P(\overline{B}_{\bf l}((\bf R_l-R_k)-(\bf R_l-M_{kl})))\\
&\leq& P(\overline{B}_{\bf l}((\bf M_{kl}-R_{k})))\\
&\leq&(m_{l_1}k_2+m_{l_2}k_1)\max\left\{ P\left( X_{\mathbf{i}}>u_{\mathbf{l},\mathbf{i}%
}\right) :\mathbf{i}\leq \mathbf{l}\right\}\\
&<\!<&\frac{m_{l_1}k_{2}}{l_1l_{2}}+\frac{m_{l_2}k_{1}}{l_{1}l_2}.
\end{eqnarray*}
Similarly, we have
$$I_3 <\!<\frac{m_{l_1}k_{2}}{l_1l_{2}}+\frac{m_{l_2}k_{1}}{l_{1}l_2}.$$
Condition $D^*(u_{\mathbf{n},\mathbf{i}})$ implies
\begin{equation*}
I_2=\left|P(B_{\bf k}({\bf R_k})\cap B_{\bf l}({\bf R_l-M_{kl}}))-P(B_{{\bf k}}({\bf R_k}))P(B_{\bf l}(\bf R_l-M_{kl}))\right|\leq \alpha_{\mathbf{l},m_{l_1},m_{l_2}}.
\end{equation*}
Noticing $\alpha_{\mathbf{l,k},m_{l_1},m_{l_2}}<\!<(\log l_1 \log l_2)^{-(\epsilon+1)}$, we obtain
$$
\left|Cov(\indi_{\left\{\bigcap_{\bf i\in R_k}\left\{X_{\bf i}\leq u_{\bf k,i}\right\}\right\}}, \indi_{\left\{\bigcap_{\bf i\in R_l-R_k}\left\{X_{\bf i}\leq u_{\bf l,i}\right\}\right\}}) \right|<\!<\alpha_{\mathbf{l},m_{l_1},m_{l_2}}+\frac{m_{l_1}k_{2}}{l_1l_{2}}+\frac{m_{l_2}k_{1}}{l_{1}l_2}.
$$
\\

\begin{lem}
Let ${\bf X}$ be a nonstationary random field such that $\left\{ n_{1}n_{2}\max\left\{ P\left( X_{\mathbf{i}}>u_{\mathbf{n},\mathbf{i}
}\right) :\mathbf{i}\leq \mathbf{n}\right\} \right\} _{\mathbf{n}\geq
\mathbf{1}}$ is bounded. Then, for $\bf k,l\in R_n$ such that $\bf k\neq \bf l$ and $u_{\bf l, i}\geq u_{\bf k, i}$,
\begin{equation*}
E\left|\indi_{\left\{\cap_{\bf i\in R_l-R_k}\left\{X_{\bf i}\leq u_{\bf l,i}\right\}\right\}}-\indi_{\left\{\cap_{\bf i\in R_l}\left\{X_{\bf i}\leq u_{\bf l,i}\right\}\right\}}\right|\leq\frac{l_1l_2-\sharp(\bf R_l-R_k)}{l_1l_2}.
\end{equation*}
\end{lem}
\textbf{Proof:}
Using the condition that $\left\{ n_{1}n_{2}\max\left\{ P\left( X_{\mathbf{i}}>u_{\mathbf{n},\mathbf{i}%
}\right) :\mathbf{i}\leq \mathbf{n}\right\} \right\} _{\mathbf{n}\geq
\mathbf{1}}$ is bounded we get
\begin{eqnarray*}
&&E\left|\indi_{\left\{\cap_{\bf i\in R_l-R_k}\left\{X_{\bf i}\leq u_{\bf l,i}\right\}\right\}}-\indi_{\left\{\cap_{\bf i\in R_l}\left\{X_{\bf i}\leq u_{\bf l,i}\right\}\right\}}\right|\\
&=&P\left(\bigcap_{{\bf i}\in {\bf R_l}-{\bf R_k}}\left\{X_{\bf i}\leq u_{\bf l,i}\right\}\right)-P\left(\bigcap_{{\bf i}\in {\bf R_l}}\left\{X_{\bf i}\leq u_{\bf l,i}\right\}\right)\\
&\leq&\sum_{\bf i\in R_l-(R_l-R_k)}P(X_{\bf i}>u_{\bf l,i})\\
&\leq& [l_1l_2-\sharp(\bf R_l-R_k)]\max\left\{ P\left( X_{\mathbf{i}}>u_{\mathbf{l},\mathbf{i}
}\right) :\mathbf{i}\leq \mathbf{l}\right\}\\
&<\!<& \frac{l_1l_2-\sharp(\bf R_l-R_k)}{l_1l_2}.
\end{eqnarray*}

The following lemma is from Tan and Wang (2014).

\begin{lem}
Let $\eta_{\bf i}$, ${\bf i}\in \mathbb{Z}_+^2$, be uniformly bounded variables. Assume that
\begin{equation*}
Var\left(\frac{1}{\log n_1 \log n_2}\sum_{\bf k \in R_n}\frac{1}{k_1k_2}\eta_{\bf k}\right)<\!<\frac{1}{(\log n_1 \log n_2)^{\epsilon+1}}.
\end{equation*}
Then
\begin{equation*}
\frac{1}{\log n_1 \log n_2}\sum_{{\bf k}\in {\bf R_n}}\frac{1}{k_1k_2}(\eta_{\bf k}-E(\eta_{\bf k}))\rightarrow 0 \ \ a.s.
\end{equation*}
\end{lem}

\textbf{Proof of Theorem 2.1:}
Let $\eta_{\bf k}=\indi_{\left\{\bigcap_{{\bf i}\leq {\bf k}}\left\{X_{\bf i}\leq u_{\bf k,i}\right\}\right\}}-E\left(\indi_{\left\{\bigcap_{{\bf i}\leq {\bf k}}\left\{X_{\bf i}\leq u_{\bf k,i}\right\}\right\}}\right)$. Then
\begin{eqnarray*}
&&Var\left(\frac{1}{\log n_1 \log n_2}\sum_{\bf k \in R_n}\frac{1}{k_1k_2}\indi_{\left\{\bigcap_{{\bf i}\leq {\bf k}}\left\{X_{\bf i}\leq u_{\bf k,i}\right\}\right\}}\right)\\
&=&E\left(\frac{1}{\log n_1 \log n_2}\sum_{\bf k \in R_n}\frac{\eta_{\bf k}}{k_1k_2}\right)^2\\
&=&\frac{1}{\log^2n_1 \log^2n_2}\left(\sum_{\bf k \in R_n}\frac{E(\eta_{\bf k}^2)}{k_1^2k_2^2}+ \sum_{{\bf k,l \in R_n},{\bf k}\neq {\bf l}}\frac{E(\eta_{\bf k}\eta_{\bf l})}{k_1k_2l_1l_2}\right)\\
&=&T_1+T_2.
\end{eqnarray*}
Since $|\eta_{\bf k}|\leq 1$, it follows that
\begin{equation*}
T_1\leq \frac{1}{\log^2n_1 \log^2n_2}\sum_{{\bf k}\in {\bf R_n}}\frac{1}{k_1^2 k_2^2}\leq \frac{K}{\log^2n_1 \log^2n_2}.
\end{equation*}
Note that for ${\bf k}\neq {\bf l}$ such that $u_{\bf k,i}<u_{\bf l,i}$,
\begin{eqnarray*}
|E(\eta_{\bf k}\eta_{\bf l})|&=&\left|Cov\left(\indi_{\left\{\bigcap_{\bf i \in R_k}\left\{X_{\bf i}\leq u_{\bf k,i}\right\}\right\}},\indi_{\left\{\bigcap_{\bf i \in R_l}\left\{X_{\bf i}\leq u_{\bf l,i}\right\}\right\}}\right)\right|\\
&\leq&\left|Cov\left(\indi_{\left\{\bigcap_{\bf i \in R_k}\left\{X_{\bf i}\leq u_{\bf k,i}\right\}\right\}},\indi_{\left\{\bigcap_{\bf i \in R_l}\left\{X_{\bf i}\leq u_{\bf l,i}\right\}\right\}}-\indi_{\left\{\bigcap_{\bf i \in R_l-R_{k}}\left\{X_{\bf i}\leq u_{\bf l,i}\right\}\right\}}\right)\right|\\
&&+\left|Cov\left(\indi_{\left\{\bigcap_{\bf i \in R_k}\left\{X_{\bf i}\leq u_{\bf k,i}\right\}\right\}},\indi_{\left\{\bigcap_{\bf i \in R_l-R_{k}}\left\{X_{\bf i}\leq u_{\bf l,i}\right\}\right\}}\right)\right|\\
&\leq&E\left|\indi_{\left\{\bigcap_{\bf i \in R_l}\left\{X_{\bf i}\leq u_{\bf l,i}\right\}\right\}}-\indi_{\left\{\bigcap_{\bf i \in R_l-R_{k}}\left\{X_{\bf i}\leq u_{\bf l,i}\right\}\right\}}\right|\\
&&+\left|Cov\left(\indi_{\left\{\bigcap_{\bf i \in R_k}\left\{X_{\bf i}\leq u_{\bf k,i}\right\}\right\}},\indi_{\left\{\bigcap_{\bf i \in R_l-R_k}\left\{X_{\bf i}\leq u_{\bf l,i}\right\}\right\}}\right)\right|.
\end{eqnarray*}
By Lemma A.2. we get
\begin{eqnarray*}
E\left|\indi_{\left\{\bigcap_{\bf i \in R_l}\left\{X_{\bf i}\leq u_{\bf l,i}\right\}\right\}}-\indi_{\left\{\bigcap_{\bf i \in R_l-R_{k}}\left\{X_{\bf i}\leq u_{\bf l,i}\right\}\right\}}\right|\leq \frac{l_1l_2-\sharp({\bf R_l-R_k})}{l_1l_2}
\end{eqnarray*}
and from Lemma A.1. we obtain
\begin{eqnarray*}
\left|Cov(\indi_{\left\{\bigcap_{\bf i \in R_k}\left\{X_{\bf i}\leq u_{\bf k,i}\right\}\right\}},\indi_{\left\{\bigcap_{\bf i \in R_l-R_k}\left\{X_{\bf i}\leq u_{\bf l,i}\right\}\right\}})\right|<\!< \alpha_{\mathbf{l},m_{l_1},m_{l_2}}+\frac{m_{l_1}k_{2}}{l_1l_{2}}+\frac{m_{l_2}k_{1}}{l_{1}l_2}.
\end{eqnarray*}
Hence
\begin{equation*}
|E(\eta_{\bf k}\eta_{\bf l})|<\!<\frac{l_1l_2-\sharp({\bf R_l-R_k})}{l_1l_2}+\alpha_{\mathbf{l},m_{l_1},m_{l_2}}+\frac{m_{l_1}k_{2}}{l_1l_{2}}+\frac{m_{l_2}k_{1}}{l_{1}l_2}.
\end{equation*}
In order to consider $T_2$, we define $A_{\bf m}=\left\{\left({\bf k},{\bf l}\right)\in {\bf R_n\times R_n}:(2m_j-1)(k_j-l_j)\geq 0, {\bf k}\neq{\bf l} \right\}$ for ${\bf m}\in\Lambda=\left\{(m_1,m_2):m_1,m_2\in \left\{0,1\right\}, {\bf m\neq 1}\right\}$. Then, we have

\begin{eqnarray*}
T_2&\leq&\frac{1}{(\log n_1 \log n_2)^2}\sum_{{{\bf m} \in \Lambda}}\sum_{({\bf k},{\bf l})\in A_{\bf m}}\frac{l_1l_2-\sharp(\bf R_l-R_k)}{l_1^2l_2^2k_1k_2}\\
&+&\frac{1}{(\log n_1 \log n_2)^2}\sum_{{{\bf m} \in \Lambda}}\sum_{({\bf k},{\bf l})\in A_{\bf m}}\frac{\alpha_{\mathbf{l},m_{l_1},m_{l_2}}+\frac{m_{l_1}k_{2}}{l_1l_{2}}+\frac{m_{l_2}k_{1}}{l_{1}l_2}}{k_1k_2l_1l_2}=:T_{21}+T_{22}.
\end{eqnarray*}
Since
\begin{eqnarray*}
T_{21}&=&\frac{1}{\log^2n_1 \log^2n_2}\underset{\underset{1\leq k_2\leq l_2\leq n_2, {\bf k\neq l}}{1\leq k_1\leq l_1\leq n_1}}%
{\sum}\bigg[\frac{k_1k_2}{l_1l_2}\times\frac{1}{k_1k_2l_1l_2}+\frac{1}{k_1k_2l_1l_2}\times\frac{k_1}{l_1}+\frac{1}{k_1k_2l_1l_2}\times\frac{k_2}{l_2}\bigg]\\
&\leq&\frac{K}{\log^2n_1 \log^2n_2}\bigg[\prod_{i=1}^2\underset{1\leq k_i\leq l_i\leq n_i}{\sum}\frac{1}{l_i^2}+\underset{1\leq k_1< l_1\leq n_1}{\sum}\frac{1}{l_1^2}\underset{1\leq l_2< k_2\leq n_2}{\sum}\frac{1}{k_2l_2}\\
&+&\underset{1\leq k_2< l_2\leq n_2}{\sum}\frac{1}{l_2^2}\underset{1\leq l_1< k_1\leq n_1}{\sum}\frac{1}{k_1l_1}\bigg]\\
&\leq&K\left(\frac{1}{\log n_1 \log n_2}+\frac{\log n_2}{\log n_1 \log n_2}+\frac{\log n_1}{\log n_1 \log n_2}\right)
\end{eqnarray*}
and
\begin{eqnarray*}
T_{22}&=&\frac{K}{(\log n_1 \log n_2)^2}\bigg[\underset{\underset{1\leq k_2\leq l_2\leq n_2, {\bf k\neq l}}{1\leq k_1\leq l_1\leq n_1}}%
{\sum}\frac{1}{k_1k_2l_1l_2(\log l_1 \log l_2)^{\epsilon_1+1}}\\
&+&\underset{1\leq k_2\leq l_2\leq n_2}{\sum}\frac{1}{k_2l_2(\log l_2)^{\epsilon_1}}\underset{1\leq l_1\leq k_1\leq n_1}{\sum}\frac{1}{k_1l_1(\log l_1)^{\epsilon_1+1}}\\
&+&\underset{1\leq k_1\leq l_1\leq n_1}{\sum}\frac{1}{k_1l_1(\log l_1)^{\epsilon_1}}\underset{1\leq l_2\leq k_2\leq n_2}{\sum}\frac{1}{k_2l_2(\log l_2)^{\epsilon_1+1}}\bigg]\\
&\leq& K(\log n_1 \log n_2)^{-(\epsilon_1+1)}
\end{eqnarray*}
we have
$$
T_2\leq K\left(\frac{1}{\log n_1 \log n_2}+\frac{\log n_2}{\log n_1 \log n_2}+\frac{\log n_1}{\log n_1 \log n_2}+\frac{1}{(\log n_1 \log n_2)^{\epsilon_1+1}}\right)
$$
and hence
$$
T_2\leq K \frac{1}{(\log n_1\log n_2)^{\epsilon+1}}, \ \ {\text {for \ some}} \ \epsilon>0.
$$
So
$$
Var\left(\frac{1}{\log n_1 \log n_2}\sum_{\bf k \in R_n}\frac{1}{k_1k_2}\indi_{\left\{\bigcap_{\bf i \in R_k}\left\{X_{\bf i}\leq u_{\bf k,i}\right\}\right\}} \right)\leq \frac{K}{(\log n_1\log n_2)^{\epsilon+1}}.
$$
The result follows by Lemma A.3. and Proposition 1.2.
\end{appendix}

\begin{appendix}\setcounter{equation}{2}

\section*{Appendix B: Proofs for Section 3}\setcounter{section}{1}

\refstepcounter{section}

The proof of Theorem 3.2 will be given through a technical lemma showing that (2) implies that
\begin{equation}
\sup_{\mathbf{1}\leq \mathbf{k}\leq \mathbf{n}}S_{\mathbf{n}}(\mathbf{R_{k}},\mathbf{R}_{\bf n}):= \sup_{\mathbf{1}\leq \mathbf{k}\leq \mathbf{n}}\underset{\mathbf{i}\in {\mathbf{R_{k}}},\mathbf{j}\in {\mathbf{R_{n}}}\atop \mathbf{i}\leq \mathbf{j}, \mathbf{i}\neq \mathbf{j}}{%
\sum }\left| r_{\mathbf{i},\mathbf{j}}\right| \exp \left( -\frac{\frac{1}{2}%
\left( u_{\mathbf{k},\mathbf{i}}^{2}+u_{\mathbf{n},\mathbf{j}}^{2}\right) }{%
1+\left| r_{{\mathbf{i},\mathbf{j}}}\right| }\right)<\!<(\log n_1 n_2)^{-(1+\epsilon)}
\end{equation}

\begin{lem}
Suppose that the covariance function $r_{{\bf i}, {\bf j}}$ satisfy $\left|r_{{\bf i}, {\bf j}}\right|<\rho_{\left|{\bf i}-{\bf j}\right|}$ for some sequence $\left\{\rho_{\bf n}\right\}_{{\bf n}\in \mathbb{N}^2-\left\{\bf 0\right\}}$ that verifies (2) for some $\epsilon>0$.
Let the constants $\{u_{\mathbf{n,i}},\mathbf{i}\leq \mathbf{n}\}_{\mathbf{n}\geq \mathbf{1}}$ be such that
$n_{1}n_{2}(1-\Phi(\lambda_{\mathbf{n}}))$ is bounded, where
$\lambda_{\mathbf{n}}=\min_{\mathbf{i}\in \mathbf{R_{n}}}u_{\mathbf{n,i}}$. Then (3) holds.
\end{lem}

We omit the proof, since it follows similar arguments to those of Lemmas 3.3-3.5 of Tan and Wang (2014).\\

\textbf{Proof of Theorem 3.2:}
We will denote the event $\left\{ X_{\mathbf{i}}\leq u_{\mathbf{n,i}%
}\right\} $ by $A_{\mathbf{i,n}}$. Using the Normal Comparison Lemma we obtain
\begin{eqnarray*}
\alpha_{\mathbf{n,k},m_{n_1},m_{n_2}}
&=&\sup_{\mathbf{1}\leq \mathbf{k}\leq \mathbf{n}}\underset{\left( \mathbf{I},\mathbf{J}\right) \in S(m_{n_{1}},m_{n_{2}})}{\sup }\left| P\left( \underset{\mathbf{i}\in \mathbf{I\wedge j}%
\in \mathbf{J}}{\bigcap }A_{\mathbf{i,k}}A_{\mathbf{j,n}}\right) -P\left(
\underset{\mathbf{i}\in \mathbf{I}}{\bigcap }A_{\mathbf{i,k}}\right) P\left(
\underset{\mathbf{j}\in \mathbf{J}}{\bigcap }A_{\mathbf{j,n}}\right) \right| \\
&\leq &\sup_{\mathbf{1}\leq \mathbf{k}\leq \mathbf{n}}\underset{\left( \mathbf{I},\mathbf{J}\right) \in S(m_{n_{1}},m_{n_{2}})}{\sup }\underset{\mathbf{i}\in {\mathbf{I}},\mathbf{j}\in {\mathbf{J}}}{%
\sum }\left| r_{\mathbf{i},\mathbf{j}}\right| \exp \left( -\frac{\frac{1}{2}%
\left( u_{\mathbf{k},\mathbf{i}}^{2}+u_{\mathbf{n},\mathbf{j}}^{2}\right) }{%
1+\left| r_{{\mathbf{i},\mathbf{j}}}\right| }\right)\\
&\leq&\sup_{\mathbf{1}\leq \mathbf{k}\leq \mathbf{n}}\underset{\left( \mathbf{I},\mathbf{J}\right)\subseteq \mathbf{R_{k}}\times \mathbf{R_{n}}}{\sup }\underset{\mathbf{i}\in {\mathbf{I}},\mathbf{j}\in {\mathbf{J}}}{%
\sum }\left| r_{\mathbf{i},\mathbf{j}}\right| \exp \left( -\frac{\frac{1}{2}%
\left( u_{\mathbf{k},\mathbf{i}}^{2}+u_{\mathbf{n},\mathbf{j}}^{2}\right) }{%
1+\left| r_{{\mathbf{i},\mathbf{j}}}\right| }\right)\\
&\leq&C\sup_{\mathbf{1}\leq \mathbf{k}\leq \mathbf{n}}\underset{\mathbf{i}\in {\mathbf{R_{k}}},\mathbf{j}\in {\mathbf{R_{n}}}\atop \mathbf{i}\leq \mathbf{j}, \mathbf{i}\neq \mathbf{j}}{%
\sum }\left| r_{\mathbf{i},\mathbf{j}}\right| \exp \left( -\frac{\frac{1}{2}%
\left( u_{\mathbf{k},\mathbf{i}}^{2}+u_{\mathbf{n},\mathbf{j}}^{2}\right) }{%
1+\left| r_{{\mathbf{i},\mathbf{j}}}\right| }\right)\\
&=&C\sup_{\mathbf{1}\leq \mathbf{k}\leq \mathbf{n}}S_{\mathbf{n}}(\mathbf{R_{k}},\mathbf{R}_{\bf n}),
\end{eqnarray*}
where $C$ is a constant. So, $D^*(u_{\mathbf{n},\mathbf{i}})$ follows from Lemma B.1. Next, we show condition $D'(u_{\mathbf{n},\mathbf{i}})$ holds.
To that end, let $\bf I \in \mathcal{E}(u_{\mathbf{n},\mathbf{i}})$. Then, we have
\begin{eqnarray*}
&&k_{n_{1}}k_{n_{2}}\underset{\mathbf{i,j}\in \mathbf{I}}{\sum }P(\overline{A%
}_{\mathbf{i,n}}\overline{A}_{\mathbf{j,n}}) \\
&\leq &k_{n_{1}}k_{n_{2}}\underset{\mathbf{i,j}\in \mathbf{I}}{\sum }\left|
P(\overline{A}_{\mathbf{i,n}}\overline{A}_{\mathbf{j,n}})-P(\overline{A}_{%
\mathbf{i,n}})P(\overline{A}_{\mathbf{j,n}})\right| +k_{n_{1}}k_{n_{2}}%
\underset{\mathbf{i,j}\in \mathbf{I}}{\sum }P(\overline{A}_{\mathbf{i,n}})P(%
\overline{A}_{\mathbf{j,n}}) \\
&\leq &k_{n_{1}}k_{n_{2}}S_{\mathbf{n}}(\mathbf{I},\mathbf{I})+k_{n_{1}}k_{n_{2}}%
\underset{\mathbf{i,j}\in \mathbf{I}}{\sum }\left( 1-\Phi (u_{\mathbf{n},%
\mathbf{i}}\right) )\left( 1-\Phi (u_{\mathbf{n},\mathbf{j}}\right) ) \\
&\leq &k_{n_{1}}k_{n_{2}}S_{\mathbf{n}}(\mathbf{R_{n}},\mathbf{R_{n}})+k_{n_{1}}k_{n_{2}}\left(
\underset{\mathbf{i}\in \mathbf{R_{n}}}{\sum }\left( 1-\Phi (u_{\mathbf{n},%
\mathbf{i}}\right) )\right) ^{2} \\
&\leq &k_{n_{1}}k_{n_{2}}S_{\mathbf{n}}(\mathbf{R_{n}},\mathbf{R_{n}})+\frac{1}{k_{n_{1}}k_{n_{2}}}%
\left(\underset{\mathbf{i}\leq \mathbf{n}}{\sum }\left( 1-\Phi (u_{\mathbf{n},%
\mathbf{i}}\right) )\right)^{2}\xrightarrow [\mathbf{n}\raini]{}0,
\end{eqnarray*}
which completes the proof of Theorem 3.2.
\bigskip

We need the following facts to prove Example 3.1, which is from Choi (2002).
The covariance function $\gamma_{n}$ satisfies the following facts
\begin{eqnarray}
\label{Tan1}
\sum_{m=0}^{n}|\gamma_{m}|^{2}\leq C n^{1-1/\log_{2}^{3}}\ \ \ \mbox{and}\ \ \
\sum_{m=0}^{n}|\gamma_{m}|^{2}\geq C \frac{n^{1-1/\log_{2}^{3}}}{\log n}
\end{eqnarray}
for some constants $C$ whose value may change form place to place. From (\ref{Tan1}) and the definition of $\gamma_{\mathbf{n}}$, it is easy to see that
\begin{eqnarray}
\label{Tan11}
\sum_{\mathbf{m}\in \mathbf{R_{n}}}|\gamma_{\mathbf{m}}|^{2}\leq C (n_{1}n_{2})^{(1-1/\log_{2}^{3})}\ \ \ \mbox{and}\ \ \
\sum_{\mathbf{m}\in \mathbf{R_{n}}}|\gamma_{\mathbf{m}}|^{2}\geq C \frac{n_{1}^{1-1/\log_{2}^{3}}}{\log n_{1}}\frac{n_{2}^{1-1/\log_{2}^{3}}}{\log n_{2}}.
\end{eqnarray}

\textbf{Proof of Example 3.1:}
We only need to show that conditions $D'\mathbb{(}u_{\mathbf{n}}\mathbb{)}$ and $D^*\mathbb{(}u_{\mathbf{n}}\mathbb{)}$ hold.
The checking of condition $D'\mathbb{(}u_{\mathbf{n}}\mathbb{)}$ is same as the proof of Theorem 3.2, so we omit it.
We will denote the event $\left\{ X_{\mathbf{i}}\leq u_{\mathbf{n}%
}\right\} $ by $B_{\mathbf{i,n}}$. Using the Normal Comparison Lemma, as for the proof Theorem 3.2, we obtain
\begin{eqnarray*}
\alpha_{\mathbf{n},m_{n_1},m_{n_2}}&=&\sup_{\mathbf{1}\leq \mathbf{k}\leq \mathbf{n}}\alpha_{\mathbf{n,k},m_{n_1},m_{n_2}}\\
&=&\sup_{\mathbf{1}\leq \mathbf{k}\leq \mathbf{n}}\underset{\left( \mathbf{I},\mathbf{J}\right) \in S(m_{n_{1}},m_{n_{2}})}{\sup }\left| P\left( \underset{\mathbf{i}\in \mathbf{I\wedge j}\in \mathbf{J}}{\bigcap }B_{\mathbf{i,k}}B_{\mathbf{j,n}}\right) -P\left(
\underset{\mathbf{i}\in \mathbf{I}}{\bigcap }B_{\mathbf{i,k}}\right) P\left(
\underset{\mathbf{j}\in \mathbf{J}}{\bigcap }B_{\mathbf{j,n}}\right) \right| \\
&\leq &C\sup_{\mathbf{1}\leq \mathbf{k}\leq \mathbf{n}}\underset{\mathbf{i}\in {\mathbf{R_{k}}},\mathbf{j}\in {\mathbf{R_{n}}}\atop \mathbf{i}\leq \mathbf{j}, \mathbf{i}\neq \mathbf{j}}{%
\sum }\left| \gamma_{\mathbf{i},\mathbf{j}}\right| \exp \left( -\frac{\frac{1}{2}\left( u_{\mathbf{k}}^{2}+u_{\mathbf{n}}^{2}\right) }{
1+\left| \gamma_{{\mathbf{i},\mathbf{j}}}\right| }\right)\\
&\leq &C\sup_{\mathbf{1}\leq \mathbf{k}\leq \mathbf{n}}k_{1}k_{2}\underset{\mathbf{0}\leq \mathbf{j}\leq \mathbf{n}, \mathbf{j}\neq \mathbf{0}}{%
\sum }\left| \gamma_{\mathbf{j}}\right| \exp \left( -\frac{\frac{1}{2}\left( u_{\mathbf{k}}^{2}+u_{\mathbf{n}}^{2}\right) }{
1+\left| \gamma_{{\mathbf{j}}}\right| }\right)=:C\sup_{\mathbf{1}\leq \mathbf{k}\leq \mathbf{n}}S_{\mathbf{n}}^{*}(\mathbf{R_{k}},\mathbf{R_{n}}).
\end{eqnarray*}


Let $\delta=\sup_{\mathbf{m}\geq \mathbf{0}, \mathbf{m}\neq\mathbf{0}}|\gamma_{\mathbf{m}}|<1$ and
$\theta_{\mathbf{n}}=\exp(\alpha u_{\mathbf{n}}^{2})$, where $\alpha$ is a constant satisfying $0<\alpha<(1-\delta)/4(1+\delta)$.
Split the term $S_{\mathbf{n}}^{*}(\mathbf{R_{k}},\mathbf{R_{n}})$ into two parts as:
$$S_{\mathbf{n}}^{*}(\mathbf{R_{k}},\mathbf{R_{n}})=\sum_{\mathbf{0}\leq \mathbf{j}\leq \mathbf{n}, \mathbf{j}\neq \mathbf{0}, \atop \chi(\mathbf{|j-i|})\leq \theta_{\mathbf{n}}}
+\sum_{\mathbf{0}\leq \mathbf{j}\leq \mathbf{n}, \mathbf{j}\neq \mathbf{0},\atop \chi(\mathbf{|j-i|})> \theta_{\mathbf{n}}}=:S_{\mathbf{n},1}^{*}+S_{\mathbf{n},2}^{*},$$
where $\chi(\mathbf{j})=\max(j_{1},1)\times\max(j_{2},1)$.
The following facts that for sufficiently large $\mathbf{n}$
\begin{eqnarray}
\label{Tan2}
\exp\left(-\frac{u_{\mathbf{n}}^{2}}{2}\right)\thicksim C\frac{u_{\mathbf{n}}}{n_{1}n_{2}}\ \ \ \mbox{and}\ \ \ u_{\mathbf{n}}\thicksim \sqrt{2\log (n_{1}n_{2})},
\end{eqnarray}
will be extensively used in the following proof. For the term $S_{\mathbf{n},1}^{*}$, using (\ref{Tan2}), we have
\begin{eqnarray*}
\sup_{\mathbf{1}\leq \mathbf{k}\leq \mathbf{n}}S_{\mathbf{n},1}^{*}
&=&\sup_{\mathbf{1}\leq \mathbf{k}\leq \mathbf{n}}k_{1}k_{2}\sum_{\mathbf{0}\leq \mathbf{j}\leq \mathbf{n}, \mathbf{j}\neq \mathbf{0}, \atop \chi(\mathbf{j})\leq \theta_{\mathbf{n}}}
\left| \gamma_{\mathbf{j}}\right| \exp \left( -\frac{\frac{1}{2}%
\left( u_{\mathbf{k}}^{2}+u_{\mathbf{n}}^{2}\right) }{%
1+\left| \gamma_{\mathbf{j}}\right| }\right)\\
&\leq&\sup_{\mathbf{1}\leq \mathbf{k}\leq \mathbf{n}}k_{1}k_{2}\sum_{\mathbf{0}\leq \mathbf{j}\leq \mathbf{n}, \mathbf{j}\neq \mathbf{0}, \atop \chi(\mathbf{j})\leq \theta_{\mathbf{n}}}
\delta \exp \left( -\frac{\frac{1}{2}\left( u_{\mathbf{k}}^{2}+u_{\mathbf{n}}^{2}\right) }{1+\delta }\right)\\
&\ll&\sup_{\mathbf{1}\leq \mathbf{k}\leq \mathbf{n}}k_{1}k_{2}\theta_{\mathbf{n}}^{2}\exp\left(-\frac{u_{\mathbf{k}}^{2}+u_{\mathbf{n}}^{2}}{2}\right)^{1/(1+\delta)}\\
&\ll&\sup_{\mathbf{1}\leq \mathbf{k}\leq \mathbf{n}}k_{1}k_{2}\theta_{\mathbf{n}}^{2}\left(\frac{u_{\mathbf{k}}}{k_{1}k_{2}}\frac{u_{\mathbf{n}}}{n_{1}n_{2}}\right)^{1/(1+\delta)}\\
&\leq&(n_{1}n_{2})^{1+4\alpha-2/(1+\delta)}(\log n_{1}n_{2})^{1/(1+\delta)}.
\end{eqnarray*}
Since $1+4\alpha-2/(1+\delta)<0$, we get $S_{1}\leq (n_{1}n_{2})^{-\kappa}$ for some $\kappa>0$.

We split the term $S_{\mathbf{n},2}^{*}$ into three parts, the first for $\mathbf{j}>\mathbf{0}$, the second for $j_{1}=0\wedge j_{2}>0$, the third for
$j_{2}=0\wedge j_{1}>0$.
We will denote them by $\mathbf{S}^{*}_{\mathbf{n},2i}$, $i=1,2,3,$ respectively.

To deal with the first case $\mathbf{j}>\mathbf{0}$, let
$$\mathbf{A}_{\mathbf{n}}=\left\{\mathbf{m}|\mathbf{1}\leq \mathbf{m}\leq \mathbf{n}, \chi(\mathbf{m})>\theta_{\mathbf{n}}, |\gamma_{\mathbf{m}}|>\frac{1}{(\log m_{1}m_{2})^{3}}\right\}.$$
 Now, we have
\begin{eqnarray*}
\sup_{\mathbf{1}\leq \mathbf{k}\leq \mathbf{n}}S_{\mathbf{n},21}^{*}&=&\sup_{\mathbf{1}\leq \mathbf{k}\leq \mathbf{n}}k_{1}k_{2}\sum_{\mathbf{j}\in \mathbf{A}_{\mathbf{n}}^{c}}\left| \gamma_{\mathbf{j}}\right| \exp \left( -\frac{\frac{1}{2}%
\left( u_{\mathbf{k}}^{2}+u_{\mathbf{n}}^{2}\right) }{%
1+\left| \gamma_{\mathbf{j}}\right| }\right)+\sup_{\mathbf{1}\leq \mathbf{k}\leq \mathbf{n}}k_{1}k_{2}\sum_{\mathbf{j}\in \mathbf{A}_{\mathbf{n}}}\left| \gamma_{\mathbf{j}}\right| \exp \left( -\frac{\frac{1}{2}%
\left( u_{\mathbf{k}}^{2}+u_{\mathbf{n}}^{2}\right) }{%
1+\left| \gamma_{\mathbf{j}}\right| }\right)\\
&=:&S_{1}+S_{2}.
\end{eqnarray*}
Since
$$\max_{\mathbf{j}\in \mathbf{A}_{\mathbf{n}}^{c}}|\gamma_{\mathbf{j}}|\leq \frac{1}{(\log \theta_{\mathbf{n}})^{3}},$$
by the same arguments as for $\mathbf{S}^{*}_{\mathbf{n},1}$, we have
\begin{eqnarray*}
S_{1}&\leq& \sup_{\mathbf{1}\leq \mathbf{k}\leq \mathbf{n}}k_{1}k_{2}n_{1}n_{2}\frac{1}{(\log \theta_{\mathbf{n}})^{3}}\exp\left(-\frac{u_{k}^{2}+u_{n}^{2}}{2(1+\frac{1}{(\log \theta_{\mathbf{n}})^{3}})}\right)\\
&\ll& \sup_{\mathbf{1}\leq \mathbf{k}\leq \mathbf{n}}k_{1}k_{2}n_{1}n_{2}\frac{1}{u_{\mathbf{n}}^{6}}\left(\frac{u_{\mathbf{k}}}{k_{1}k_{2}}\frac{u_{\mathbf{n}}}{n_{1}n_{2}}\right)^{1+\frac{1}{\alpha^{3}u_{\mathbf{n}}^{6}}}\\
&\ll& (n_{1}n_{2})^{-\frac{2}{\alpha^{3}u_{\mathbf{n}}^{6}}}(u_{\mathbf{n}})^{-4+\frac{2}{\alpha^{3}u_{\mathbf{n}}^{6}}}\\
&\ll& (\log n_{1}n_{2})^{-2}.
\end{eqnarray*}

Now we consider the term $S_{2}$. Let $\beta=1-1/(2\log_{2}^{3})$.
Form the definition of $\gamma_{\mathbf{m}}$, we have
$$\delta':=\sup_{\mathbf{m}\in \mathbf{A}_{\mathbf{n}}}|\gamma_{\mathbf{m}}|\leq \sup_{\mathbf{m}\in \mathbf{A}_{\mathbf{n}}}\left(\frac{1}{\log m_{1}\log m_{2}}\right)^{1/2}
\leq \sup_{\mathbf{m}\in \mathbf{A}_{\mathbf{n}}}\left(\frac{1}{\log m_{1}m_{2}}\right)^{1/2}\leq \left(\frac{1}{\log \theta_{\mathbf{n}}}\right)^{1/2}$$
As in Choi (2002), we clam that $card(\mathbf{A}_{\mathbf{n}})=O((n_{1}n_{2})^{\beta})$.  If not, $|\gamma_{\mathbf{m}}|>\frac{1}{(\log m_{1}m_{2})^{3}}$ on a set of size $O((n_{1}n_{2})^{\beta})$ and thus
$$\sum_{\mathbf{m}\in \mathbf{R_{n}}}|\gamma_{\mathbf{m}}|^{2}\geq \sum_{\mathbf{m}\in \mathbf{A}_{\mathbf{n}}}|\gamma_{\mathbf{m}}|^{2}\geq C\frac{(n_{1}n_{2})^{\beta}}{(\log n_{1}n_{2})^{6}}$$
contradicting (\ref{Tan11}).
Hence
\begin{eqnarray*}
S_{2}&=&\sup_{\mathbf{1}\leq \mathbf{k}\leq \mathbf{n}}k_{1}k_{2}\sum_{\mathbf{j}\in \mathbf{A}_{\mathbf{n}}}\left| \gamma_{\mathbf{j}}\right| \exp \left( -\frac{\frac{1}{2}%
\left( u_{\mathbf{k}}^{2}+u_{\mathbf{n}}^{2}\right) }{%
1+\left| \gamma_{\mathbf{j}}\right| }\right)\\
&\leq&\sup_{\mathbf{1}\leq \mathbf{k}\leq \mathbf{n}}k_{1}k_{2}(n_{1}n_{2})^{\beta}\frac{1}{(\log \theta_{\mathbf{n}})^{1/2}}\exp\left(-\frac{u_{\mathbf{k}}^{2}+u_{\mathbf{n}}^{2}}{2(1+\delta')}\right)\\
&\ll& (n_{1}n_{2})^{1+\beta-\frac{2}{1+\delta'}}(u_{\mathbf{n}})^{\frac{2}{1+\delta'}-1}\\
&\ll& (n_{1}n_{2})^{2-\frac{1}{2\log_{2}^{3}}-\frac{2}{1+\delta'}}(u_{\mathbf{n}})^{\frac{2}{1+\delta'}-1}\\
&\ll& (n_{1}n_{2})^{-\varepsilon},
\end{eqnarray*}
for some $\varepsilon>0$.

Next, we deal with the second case $j_{1}=0\wedge j_{2}>0$. If $n_{2}\leq\theta_{\mathbf{n}}$, by the same argument as for $\mathbf{S}^{*}_{\mathbf{n},1}$,
we can show
$$\sup_{\mathbf{1}\leq \mathbf{k}\leq \mathbf{n}}S_{\mathbf{n},22}^{*}\leq (n_{1}n_{2})^{-\varepsilon}$$
for some $\varepsilon>0$.
If $n_{2}>\theta_{\mathbf{n}}$, let
$$\mathbf{B}_{\mathbf{n}}=\left\{(0,m_{2})|1\leq m_{2}\leq n_{2}, m_{2}>\theta_{\mathbf{n}}, |\gamma_{(0,m_{2})}|>\frac{1}{(\log m_{2})^{3}}\right\}.$$

 Now, we have
\begin{eqnarray*}
\sup_{\mathbf{1}\leq \mathbf{k}\leq \mathbf{n}}S_{\mathbf{n},22}^{*}&=&\sup_{\mathbf{1}\leq \mathbf{k}\leq \mathbf{n}}k_{1}k_{2}\sum_{\mathbf{j}\in \mathbf{B}_{\mathbf{n}}^{c}}\left| \gamma_{\mathbf{j}}\right| \exp \left( -\frac{\frac{1}{2}%
\left( u_{\mathbf{k}}^{2}+u_{\mathbf{n}}^{2}\right) }{%
1+\left| \gamma_{\mathbf{j}}\right| }\right)+\sup_{\mathbf{1}\leq \mathbf{k}\leq \mathbf{n}}k_{1}k_{2}\sum_{\mathbf{j}\in \mathbf{B}_{\mathbf{n}}}\left| \gamma_{\mathbf{j}}\right| \exp \left( -\frac{\frac{1}{2}%
\left( u_{\mathbf{k}}^{2}+u_{\mathbf{n}}^{2}\right) }{%
1+\left| \gamma_{\mathbf{j}}\right| }\right)\\
&=:&S_{3}+S_{4}.
\end{eqnarray*}

Since
$$\max_{\mathbf{j}\in \mathbf{B}_{\mathbf{n}}^{c}}|\gamma_{\mathbf{j}}|\leq \frac{1}{(\log \theta_{\mathbf{n}})^{3}},$$
by the same arguments as for $S_{1}$, we have
\begin{eqnarray*}
S_{3}&\leq& \sup_{\mathbf{1}\leq \mathbf{k}\leq \mathbf{n}}k_{1}k_{2}n_{2}\frac{1}{(\log \theta_{\mathbf{n}})^{3}}\exp\left(-\frac{u_{k}^{2}+u_{n}^{2}}{2(1+\frac{1}{(\log \theta_{\mathbf{n}})^{3}})}\right)\\
&<& \sup_{\mathbf{1}\leq \mathbf{k}\leq \mathbf{n}}k_{1}k_{2}n_{1}n_{2}\frac{1}{u_{\mathbf{n}}^{6}}\left(\frac{u_{\mathbf{k}}}{k_{1}k_{2}}\frac{u_{\mathbf{n}}}{n_{1}n_{2}}\right)^{1+\frac{1}{\alpha^{3}u_{\mathbf{n}}^{6}}}\\
&\ll& (n_{1}n_{2})^{-\frac{2}{\alpha^{3}u_{\mathbf{n}}^{6}}}(u_{\mathbf{n}})^{-4+\frac{2}{\alpha^{3}u_{\mathbf{n}}^{6}}}\\
&\ll& (\log n_{1}n_{2})^{-2}.
\end{eqnarray*}

Now we consider the term $S_{4}$.
Noting that $\gamma_{\mathbf{m}}=\gamma_{m_{1}}\gamma_{m_{2}}$ and $\gamma_{0}=1$, we have
$$\delta'':=\sup_{\mathbf{m}\in \mathbf{B}_{\mathbf{n}}}|\gamma_{\mathbf{m}}|\leq \sup_{\mathbf{m}\in \mathbf{B}_{\mathbf{n}}}\left(\frac{1}{\log m_{2}}\right)^{1/2}\leq \left(\frac{1}{\log \theta_{\mathbf{n}}}\right)^{1/2}$$
As in Choi (2002), we clam that $card(\mathbf{B}_{\mathbf{n}})=O((n_{2})^{\beta})$.  If not, $|\gamma_{\mathbf{m}}|>\frac{1}{(\log m_{2})^{3}}$ on a set of size $O((n_{2})^{\beta})$ and thus
$$\sum_{m_{2}=1}^{n_{2}}|\gamma_{m_{2}}|^{2}=\sum_{m_{2}=1}^{n_{2}}|\gamma_{(0,m_{2})}|^{2}\geq \sum_{\mathbf{m}\in \mathbf{B}_{\mathbf{n}}}|\gamma_{\mathbf{m}}|^{2}\geq C\frac{(n_{2})^{\beta}}{(\log n_{2})^{6}}$$
contradicting (\ref{Tan1}).
Hence
\begin{eqnarray*}
S_{4}&=&\sup_{\mathbf{1}\leq \mathbf{k}\leq \mathbf{n}}k_{1}k_{2}\sum_{\mathbf{j}\in \mathbf{B}_{\mathbf{n}}}\left| \gamma_{\mathbf{j}}\right| \exp \left( -\frac{\frac{1}{2}\left( u_{\mathbf{k}}^{2}+u_{\mathbf{n}}^{2}\right) }{%
1+\left| \gamma_{\mathbf{j}}\right| }\right)\\
&\leq&\sup_{\mathbf{1}\leq \mathbf{k}\leq \mathbf{n}}k_{1}k_{2}(n_{2})^{\beta}\frac{1}{(\log \theta_{\mathbf{n}})^{1/2}}\exp\left(-\frac{u_{\mathbf{k}}^{2}+u_{\mathbf{n}}^{2}}{2(1+\delta')}\right)\\
&<& (n_{1}n_{2})^{1+\beta-\frac{2}{1+\delta'}}(u_{\mathbf{n}})^{\frac{2}{1+\delta'}-1}\\
&\ll& (n_{1}n_{2})^{2-\frac{1}{2\log_{2}^{3}}-\frac{2}{1+\delta'}}(u_{\mathbf{n}})^{\frac{2}{1+\delta'}-1}\\
&\ll& (n_{1}n_{2})^{-\varepsilon},
\end{eqnarray*}
for some $\varepsilon>0$.
Likewise we can bound the third case $j_{2}=0\wedge j_{1}>0$.
 Thus condition $D^*\mathbb{(}u_{\mathbf{n}}\mathbb{)}$ holds.
\end{appendix}

\bigskip
{\bf Acknowledgement}: The authors would like to thank the
referee  for several corrections and important suggestions which
significantly improved this paper. Lu\'{i}sa's work was supported by National Foundation of Science and Technology through UID/MAT/00212/2013.
Tan's work was supported by National Science Foundation of China
(No. 11501250), Natural Science Foundation of Zhejiang Province of China (No. LQ14A010012, LY15A010019
).


\begin{thebibliography}{30}

\bibitem{Berkes} {\small Berkes, I., Cs\'{a}ki, E. (2001). A universal result in almost sure central limit theorem. \textit{Stoch. Process. Appl.}, 94, 105-134.}

\bibitem{Brosamler} {\small Brosamler, G.A. (1988). An almost sure everywhere central limit theorem. \textit{Math. Proc. Camb. Phil. Soc.}, 104, 561-574.}

\bibitem{Chen} {\small Chen, S., Lin, Z. (2006). Almost sure max-limits for nonstationary Gaussian sequence. \textit{Statist. Probab. Lett.}, 76, 1175-1184.}

\bibitem{Cheng} {\small Cheng, S., Peng, L., Qi, Y. (1998). Almost sure convergence in extreme value theory. \textit{Math. Nachr.}, 190, 43-50.}

\bibitem{} Choi, H., (2002). Central limit theory and extemes of random fields. Phd Dissertation
in Univ. of North Carolina at Chapel Hill.

\bibitem{Choi}  {\small Choi, H. (2010). Almost sure limit theorem for stationary Gaussian random fields. \textit{Journal of the Korean Statistical Society}, 39, 475-482.}

\bibitem{Csaki}  {\small Cs\'{a}ki, E., Gonchigdanzan, K. (2002). Almost sure limit theorems for the maximum of stationary Gaussian sequences. \textit{Statist. Probab. Lett.}, 58, 195-203.}


\bibitem{Fahrner}  {\small Fahrner, I., Stadm\"{u}ller, U. (1998). On almost sure max-limit theorems. \textit{Statist. Probab. Lett.}, 37, 229-236.}


\bibitem{Hus1}  {\small H\"{u}sler, J. (1986) Extremes of nonstationary random sequences. \emph{Journal of Applied Probability}, 23, 937-950. }

\bibitem{Lacey}  {\small Lacey, M.T., Philipp, W. (1990). A note on the almost sure central limit theorem. \textit{Statist. Probab. Lett.}, 9, 201-205.}

\bibitem{Leadbetter2}  {\small Leadbetter and Rootz\'{e}n, H. (1998) On
extreme values in stationary random fields. \emph{Stochastic processes and
related topics}, 275-285, Trends Math. Birkhauser Boston, Boston. }

\bibitem{Pereira}  {\small Pereira, L. and Ferreira, H. (2005) Extremes of
quasi-independent random fields and clustering of high values. \emph{%
Proceedings of 8th WSEAS International Conference on Applied Mathematics},104-109. }

\bibitem{Pereira2}  {\small Pereira, L. and Ferreira, H. (2006) Limiting
crossing probabilities of random fields. \emph{Journal of Applied Probability%
}}, 3, 884-891.


\bibitem{Schatte}  {\small Schatte, P. (1988) On strong versions of the central limit theorem. \textit{Math. Nachr.}, 137, 249-256.}

\bibitem{Zuoxiang}  {\small Peng, Z., Nadarajah, S. (2011) Almost sure limit theorems for Gaussian sequences. \textit{Theory Probab. Appl.}, 55, 361-367.}



\bibitem{Tan14} {\small Tan, Z., Wang, Y. (2014) Almost sure asymptotics for extremes of non-stationary Gaussian random fields. \textit{Chinese Annals of Mathematics, Series B}, 35, 125-138.}

\end{thebibliography}
\end{document}